\documentclass[12pt]{amsart}

\newcommand{\Z}{{\mathbb Z}}
\newcommand{\F}{{\mathbb F}}
\DeclareMathOperator{\Real}{Re}

\theoremstyle{remark}
\newtheorem* {remark}     {Remark}

\begin{document}

\title{On the genera of $X_0(N)$}

\author[Csirik]{J\'anos A.\ Csirik}
\address{AT\&T Shannon Lab,
P.O.\ Box 971,
180 Park Ave,
Florham Park, NJ 07932--0971}
\email{janos@research.att.com}

\author[Wetherell]{Joseph L.\ Wetherell}
\address{Department of Mathematics,
University of Southern California,
Los Angeles, CA 90089--1113}
\email{jlwether@alum.mit.edu}

\author[Zieve]{Michael E.\ Zieve}
\address{Center for Communications Research,
         29 Thanet Rd.,
         Princeton, NJ 08540--3699}
\email{zieve@idaccr.org}

\thanks{The first author thanks A.\thinspace{}M.\ Odlyzko for explaining how 
to use Dirichlet series to calculate the cumulative average of 
$\mu/N$ (Section~\ref{sec-average}), and
also thanks H. Zhu for bringing this topic to his attention.
The third author thanks A. Granville, C. Pomerance, and J. Vanderkam
for helpful correspondence.  The second and third authors were supported in part
by NSF Mathematical Sciences Postdoctoral Research Fellowships.
The authors used the computer package {\tt GP/PARI}
for various computations related to this note.}

\begin{abstract}
Let $g_0(N)$ be the genus of the modular curve $X_0(N)$.  We record
several properties of the sequence $\{g_0(N)\}$.
Even though the average size of $g_0(N)$ is $(1.25/\pi^2)N$,
a random positive integer has probability zero of being
a value of $g_0(N)$.
% (i.e., the set of integers of the form  $g_0(N)$ has density zero)
Also, if $N$ is a random positive integer then $g_0(N)$ is odd with 
probability one.
\end{abstract}

\date{August 26, 2000}

\maketitle

%####################################################################
%####################################################################
%####################################################################

\section{Introduction}

For each positive integer $N$, the modular curve $X_0(N)$ parametrizes
elliptic curves together with a cyclic subgroup of order $N$ (for more
details, see~\cite{Ro} and \cite{Sh}).  
% Note: the moduli interpretation is Prop.1 (p.46) of Rohrlich.
% In fact, Rohrlich's Section 1.3 is titled ``Moduli interpretation
% of X_0(N) and X_1(N)''
The genus $g_0(N)$ of $X_0(N)$ tends to infinity
as the level $N$ increases.  We will examine the sequence $\{g_0(N)\}$
more closely.  Among other things, we will show that the average size of
$g_0(N)$ is $1.25N/\pi^2$; however, the values
of $g_0(N)$ form a density-zero subset of the integers.
These two properties imply that there is much collapsing under the
map $N\mapsto g_0(N)$: for instance, there are integers whose preimage
under this map is arbitrarily large.

The results of this note are:
\begin{enumerate}
\item Upper and lower bounds, including the asymptotic results \newline
$0<\liminf g_0(N)/N < \limsup g_0(N)/(N\log\log N)<\infty$.
\item Average behavior: $\lim_{B\to\infty} (1/B)\sum_{N\leq B} g_0(N)/N = 
1.25/\pi^2$.
\item Natural density:  $\{g_0(N)\}$ is a density zero subset of the integers.
\item Non-uniformity of $g_0(N)$ modulo a fixed prime $p$: for instance, $g_0(N)$ is
odd with probability 1, and (for a fixed odd prime $p$) the probability that
$g_0(N)\equiv1\pmod{p}$ is much less than $1/p$.
\end{enumerate}

This note arose from a question about curves over finite fields.
Typically one estimates the number of points on such a curve by
means of the Weil (upper) bound.  However, in various situations this bound
can be improved, which leads to the converse question of producing
curves with many points (of a given genus, over a given field)
in order to see how far the Weil bound can be improved.
It is known that, if the prime $p$ does not divide $N$, then
$X_0(N)$ has many rational points over the finite field $\F_{p^2}$.
Covers of these curves have been used to show that there are
curves with many points (over $\F_{p^2}$) in every genus~\cite{EHKPWZ}.
It is natural to ask whether covers are needed to prove the latter
result, or whether the genera of $X_0(N)$ already achieve all sufficiently
large integer values.  The present note grew from our proof that
infinitely many positive integers do not occur in the sequence $\{g_0(N)\}$.

% More precisely: if p doesn't divide N, then for k=GF(p^2),  #X_0(N)(k)  is
% g*(p-1) + o_p(g),  where g=genus of X_0(N).
% For p fixed and g growing, the main term in this estimate is the largest
% posible for any curve of that genus, according to the Drinfeld-Vladut bound.
% So, if the curves X_0(N) (with N coprime to p) attain every sufficiently large genus,
% then we would know that  lim_{g->infinity} N_{p^2}(g)/g = p-1.
% But these curves X_0(N) do not attain every large genus, so we cannot
% draw this conclusion.  Instead, the paper \cite{EKPWZ} shows that double covers
% of the curves X_0(N) (with N coprime to p) achieve every genus, and as a 
% consequence one has  liminf_g N_{p^2}(g)/g >= (p-1)/2.

The authors suspect that some of the results of this note have been discovered
previously, but apparently they have not been published.\footnote{Indeed, we
recently saw an unpublished manuscript of S. Wong which contains some of the
results of this note.}

% More precisely, after seeing a version of this note in the xxx archives, 
% Siman Wong sent us his 1995 manuscript which proves a lower bound on g_0(N)
% (namely, g_0(N) > N/12 - (3/2)*\sqrt{N}) and shows that the set of integers
% N for which g_0(N) is NOT 1 mod 2^m  has density zero.

%####################################################################
%####################################################################

\section{The genus of $X_0(N)$}
\label{sec-def}

For a positive integer $N$, let $g_0(N)$ denote the genus of $X_0(N)$.
By \cite[Prop.~1.43]{Sh}, we have 
$$
g_0(N) = 1 + \frac{\mu}{12}-\frac{\nu_2}{4}-\frac{\nu_3}{3}-\frac{\nu_\infty}{2}
$$
where $\nu_2$ and $\nu_3$ are the numbers of solutions in $\Z/N\Z$ of the
equations $x^2+1=0$ and $x^2+x+1=0$, respectively; 
$\nu_\infty=\sum_d \varphi(\text{gcd}(d,N/d))$; and $\mu=N\sum_a 1/a$.
Here $\varphi$ denotes Euler's totient function, $d$ ranges over the positive 
divisors of $N$, and $a$ ranges over the squarefree positive divisors of $N$.  
Likewise, let $p$ range over the prime divisors of $N$, and write $N=\prod_p p^{r_p}$.
Then $\mu = \prod_p (p+1)p^{r_p-1}$.  Also $\nu_2$ is 0 if $4\mid N$
or if some $p\equiv 3\pmod{4}$, and otherwise $\nu_2=2^s$ where
$s$ is the number of $p\equiv 1\pmod{4}$; similarly, $\nu_3$ is 0 if $9\mid N$
or if some $p\equiv 2\pmod{3}$, and otherwise $\nu_3=2^t$ where $t$ is
the number of $p\equiv 1\pmod{3}$.  Finally, $\nu_\infty=\prod_p \theta(p,r_p)$,
where we define $\theta(p,2R+1) = 2p^R$ and $\theta(p,2R) = (p+1)p^{R-1}$.

%####################################################################
%####################################################################

\section{Bounds on $g_0(N)$}
\label{sec-bounds}

We now give general upper and lower bounds on $g_0(N)$.
For the lower bound we show that
$$ g_0(N) \ge (N-5\sqrt{N}-8)/12,$$
with equality if and only if $N=p^2$ where $p$ is a prime
congruent to 1~(mod~12).  For the upper bound we show
the asymptotic result
$$ \limsup_{N\to\infty} \frac{g_0(N)}{N\log\log N} =
   \frac{e^{\gamma}}{2\pi^2}$$ 
and the explicit bound
$$ g_0(N) < N\frac{e^\gamma}{2\pi^2}(\log\log N + 2/\log\log N)
\quad\text{for $N>2$},$$
where $\gamma=0.5772...$ is Euler's constant.
% Definition: $\gamma=\lim_{n\to\infty}(1+1/2+...+1/n-\log{n})$.
% Note that $2\pi^2/12 = \zeta(2)$.  

We start by proving the lower bound $g_0(N)\geq (N-5\sqrt{N}-8)/12$.
If $N$ is a prime power (or $N=1$), it is easy to prove this
bound and to show that equality occurs precisely when $N$ is the
square of a prime congruent to 1~(mod~12).
% Proof when N is a prime power N=p^r:
% if r=1 then  N >= 1 + (p+1)/12 - 2/4 - 2/3 - 2/2 = (p-13)/12 > (p-5sqrt{p}-8)/12.
% if r=2 then  N >= 1 + (p^2+p)/12 - 2/4 - 2/3 - (p+1)/2 = (p^2-5p-8)/12,
%    with equality if and only if  p=1 (mod 12).
% if r>2 then  N >= 1 + mu/12 - 2/4 - 2/3 - mu/(2*sqrt{N})
%                 = (1/12)*(p^r + p^(r-1) - 6 p^(r/2) - 6 p^(r/2-1) - 2)
%                 = (1/12)*(p^r - 5p^(r/2) - 8 +  (p^(r/2-1) - 1)*(p^(r/2) - 6))
%        so the desired inequality is true whenever  N > 36.  
% Finally, we directly verify the result in case N=1, 8, 16, 27, 32.
It is also easy to check the lower bound for $N$ less than 2000, so
we may assume in the sequel that $N \ge 2000$.

Suppose that $N$ has exactly two distinct prime factors $p$ and
$q$.  Then $\mu - N > N/p + N/q  > 2\sqrt{N}$.  Since
$\nu_\infty \leq \mu/\sqrt{N}$ and $\nu_2, \nu_3 \le 4$, we have 
$g_0(N)  \ge 1 + \mu/12 - 4/4 - 4/3 - \mu/(2\sqrt{N}).$
Thus,
\begin{align*}
  12g_0(N) 
          & \ge \mu (1 - 6/\sqrt{N}) - 16 \\
          & > (N+2\sqrt{N}) (1 - 6/\sqrt{N}) - 16 \\
%         & = N - 4\sqrt{N} - 28 \\
          & > N - 5\sqrt{N} - 8.  %\qquad \text{(since $N > 400$).}
\end{align*}

Finally, consider the case where $N$ has at least 3 distinct prime
factors.  In this case we observe that $\mu - N > 3N^{2/3}$ and
$\nu_2, \nu_3 \le \nu_\infty$.  It follows that $g_0(N) \ge 1 +
\mu/12 - (1/4+1/3+1/2)\nu_\infty$; hence,
\begin{align*}
  12g_0(N) 
     &\ge 12 + \mu - 13\mu/\sqrt{N} \\
     &> (N-5\sqrt{N}-8) + (\mu - N)(1-13/\sqrt{N}) - 8\sqrt{N} + 20\\
     &> (N-5\sqrt{N}-8) + 3N^{4/6} - 39N^{1/6} - 8N^{3/6} + 20.
\end{align*} 
The largest real root of $3x^4 - 8x^3 - 39x + 20$ is
$3.5495... < 2000^{1/6}$.  The lower bound on $g_0(N)$ follows.

We now prove the limsup result.  We have already noted that $\nu_2$,
$\nu_3$, and $\nu_\infty$ are each bounded by $\mu/\sqrt{N}$; thus,
it suffices to show that
%the desired result is equivalent to
$$\limsup_{N\to\infty}\frac{\mu}{N\log\log N}=\frac{6}{\pi^2}e^\gamma.$$ 
For any $x\geq 2$ let $N_x = \prod_{p\leq x} p$.
Among all $N$ with $N_x\leq N<N_{x+1}$ (for fixed $x$), the case
$N=N_x$ maximizes $\mu/N$ and minimizes $\log\log N$.  
Thus, we need only consider values $N=N_x$.
For these $N$, Mertens showed in 1874 that
$\mu/N$ is asymptotic to $(\log x)e^\gamma 6/\pi^2$ \cite[p.~429]{HW};
%
% Actually Mertens' theorem says that  
% $$\prod_{p\leq x} (1-1/p)=(e^{-\gamma})/\log x + o(1/\log x);$$
% since $\prod_p (1-1/p^2)=6/\pi^2$ (Euler, 1755), it follows that
% $\prod_{p\leq x} (1+1/p)$ is asymptotic to $(\log x)e^\gamma 6/\pi^2$.
%
the limsup result follows, since
$\log\log N_x$ is asymptotic to $\log x$ \cite[p.~341]{HW}.
%
% The citation is to Chebyshev's 1854 result that there are positive 
% constants c,d with  c*x < SUM_{p<x} log(p) < d*x.  
% Actually Chebyshev's result gives vastly more than we need.
%

The upper bound on $g_0(N)$ is easily verified for $N<210$, so we now
assume $N\geq 210$.  Note that $g_0(N)\leq \mu/12$, so our upper bound
follows from the inequality
$\mu/N < (\log\log N + 2/\log\log N)e^\gamma/(2\pi^2)$.
As above, it suffices to prove this inequality when $N=N_x$ (since $N\geq 
210=N_7$).
But in this case the inequality follows easily from results in~\cite{RS}.

\section{Data for small $N$}
\label{sec-data}

We now determine the first few positive integers $n$ which do not occur
as $g_0(N)$ for any $N$.  The lower bound of the previous section implies that,
if $n=g_0(N)$, then $N < 12n+18\sqrt{n}+40$.  
%
% Proof: we know  12n >= N - 5\sqrt{N} - 8,  so
%          (\sqrt{N} - 5/2)^2 <= 12n + 57/4
%            \sqrt{N} <= 5/2 + \sqrt{12n + 57/4}
%             N <= 25/4 + 5\sqrt{12n + 57/4} + 12n + 57/4
%                = 12n + 41/2 + 5\sqrt{12n + 57/4}
%   and this last is  < 12n + 18\sqrt{n} + 40   if and only if
%       5\sqrt{12n + 57/4} < 18\sqrt{n} + 39/2
%       25*(12n + 57/4) < 324n + 702\sqrt{n} + 1521/4
%                    0 < 24n + 702\sqrt{n} + 24
%   which is visibly true.  
%
So, to determine whether $n$
occurs as $g_0(N)$, we just check all levels $N$ up to this bound.
Doing this by computer, we find the first few missed values (i.e., positive
integers $n$ not of the form $g_0(N)$):
150, 180, 210, 286, 304, 312, ...  
Note that all of these are even; in fact, the first several thousand missed values
are even.  It is easy to 
describe all positive integers $N$ for which $g_0(N)$ is even.
They are given in the following list, where $p$ denotes a prime and $r$ denotes
a positive integer:
\begin{enumerate}
\item $N=1,2,3,4,8$ or $16$;
\item $N=p^r$ where $p\equiv5\pmod 8$;
\item $N=p^r$ where $p\equiv7\pmod 8$ and $r$ odd;
\item $N=p^r$ where $p\equiv3\pmod 8$ and $r$ even;
\item $N=2p^r$ where $p\equiv\pm3\pmod 8$;
\item $N=4p^r$ where $p\equiv3\pmod 4$ and $r$ odd.
\end{enumerate}
It follows that, if $N$ is a randomly chosen positive integer, then
$g_0(N)$ is odd with probability 1.  This shows that a randomly
chosen even positive integer has probability zero of occurring as
a value of $g_0(N)$.

Looking further in the list of integers not of the form $g_0(N)$, 
we do eventually find some odd values, the first one occuring at the 
$3885$th position.  There are four
such up to $10^5$ (out of 9035 total missed values), namely 49267, 
74135, 94091, 96463.  
In Section~\ref{sec-density}
we will show that the paucity of odd missed values is not a
general phenomenon, but instead an accident caused by the fact that
`small' numbers do not have enough prime factors.
In particular, we will see that there are infinitely many positive odd
integers not of the form $g_0(N)$, and in fact the set of integers of the
form $g_0(N)$ has density zero in the set of all nonnegative integers.

%####################################################################
%####################################################################

\section{Average size of $g_0(N)$}
\label{sec-average}

We have shown that $g_0(N)$ is sometimes as small as $N/12$, and
sometimes as large as $cN\log\log N$.  We now determine the usual behavior.
More precisely, we show that
$$\frac1B\sum_{N=1}^B g_0(N) = \frac{5}{8\pi^2}B+o(B).$$
By Abel's lemma~\cite[p.~65]{Se}, this result is equivalent to the following:
\begin{equation}
\label{avg}
\lim_{B\to\infty} \frac1B\sum_{N=1}^B \frac{g_0(N)}N = \frac{5}{4\pi^2}=0.12665...
\end{equation}
In Section~\ref{sec-bounds}, we showed that
$g_0(N)/N=(1/12)\sum_{a\mid N}(1/a) + o(1)$,
and we can ignore the error term since it contributes nothing to the limit.
Here we use the convention that $p$ is always prime,  and $a$ is 
always a squarefree positive integer.
The following computation proves (\ref{avg}):
\begin{eqnarray*}
\lim_{B\to\infty} (1/B)\sum_{N=1}^B \sum_{a\mid N} 1/a
%&= \lim_{B\to\infty} (1/B)\sum_{a\leq B} (1/a)\lfloor B/a\rfloor \\
= \lim_{B\to\infty} (1/B)\sum_{a\leq B} (B/a)(1/a) \\
%&= \lim_{B\to\infty} \sum_{a\leq B} 1/a^2 \\
= \sum_{a<\infty} 1/a^2 
= \zeta(2) / \zeta(4) 
%&= (\pi^2/6)  /  (\pi^4/90) \\
= 15/\pi^2.
\end{eqnarray*}

A.\thinspace{}M.\ Odlyzko showed us the following alternative proof of (\ref{avg}):
for any complex $s$ with $\Real(s) > 1$, the series 
$F(s):=\sum_{N=1}^\infty (\mu/N)N^{-s}$ converges,
% (it converges absolutely, due to the bounds of Section 3)
and we have
\begin{align*}
F(s)&=\sum_{N=1}^\infty N^{-s}\prod_{p\mid N} (1+1/p)
= \prod_p \left(1+(1+1/p)(p^{-s}+p^{-2s}+\cdots)\right) \\
%&= \prod_p(1+p^{-s}(1+1/p)/(1-p^{-s}) ) \\
%&= \zeta(s)\prod_p (1-1/p^s+1/p^s+1/p^{s+1}) \\
&= \zeta(s)\prod_p (1+1/p^{s+1}) 
= \zeta(s)\zeta(s+1)/\zeta(2s+2).
\end{align*}
It follows that $F(s)$ has a simple pole at $s=1$ with residue $\zeta(2)/\zeta(4)$,
and that $F(s)$ is regular on $\Real(s)=1$ except at $s=1$.  
Now (\ref{avg})
follows at once from Ikehara's Tauberian theorem~\cite[p.~311]{La}.

%####################################################################
%####################################################################

\section{Density of $\{g_0(N)\}$}
\label{sec-density}

In this section we consider the (natural) density of the set
$\{g_0(N)\}$ as a subset of the non-negative integers.  We have
already seen that $g_0(N)$ is almost never even, so this density (if
it exists) is at most $1/2$.  In fact, we will show that the set
$\{g_0(N)\}$ has density zero---in other words, we will show that
$\lim_{x\to\infty} (\#S(x))/x=0$, where
$$
S(x)=\{n\in\Z\colon n\leq x \text{ and } n=g_0(N)\text{ for some }N\}.
$$

Let $N$ be a positive integer and suppose that $N$ is divisible by at
least $s>2$ distinct odd primes.  The formulas in Section~\ref{sec-def}
imply that $2^{s-1}$ divides $\nu_3$ and $2^s$ divides each of 
$\mu$, $\nu_2$, and $\nu_{\infty}$.  It follows that 
$g_0(N) \equiv 1 \pmod{2^{s-2}}$.

We now show that $\#S(x)\leq x/2^d + o_d(x)$ for each positive integer~$d$; 
this implies that $\#S(x)=o(x)$, as desired.
Fix a positive integer $d$.  Clearly the number of $n\in S(x)$ with
$n\equiv 1\pmod{2^d}$ is less than $x/2^d+1$.  It remains to show
that the number of $n\in S(x)$ with $n\not\equiv 1\pmod{2^d}$ is $o(x)$.
Each such $n$ has the form $g_0(N)$ where $N$ has at most $d+1$ distinct
odd prime divisors and $N<12x+18\sqrt{x}+40$.
And the number of such $N$ (hence also the number of such $n$) is well-known
to be $o(x)$ \cite[p.~356]{HW}.

\begin{remark}
We have shown that $\#S(x)=o(x)$.  It would be interesting to obtain a
more precise estimate for $\#S(x)$.  We make some preliminary observations
in this vein.  First, the number of even integers in $S(x)$ is asymptotic 
to $9x/\log x$.
(This follows from the classification of $N$ for which $g_0(N)$ is even,
together with the fact that, for any fixed nonzero integers $a,b$, there
are only $o(x/\log x)$ primes $p<x$ for which $|ap+b|$ is 
prime \cite[Cor.~2.4.1]{HRi}.)
%
% (The details are given above, immediately following the description
%  of the N for which g_0(N) is even.)
%
We can prove the following upper and lower bounds on $\#S(x)$ for $x\geq 3$:
$$ \frac{x}{\log x} e^{a(\log\log\log x)^2} \ll \#S(x) \ll \frac{x}
  {(\log x)^b(\log\log x)^c}.$$
(For positive $f(x),g(x)$, the notation $f(x)\ll g(x)$ is equivalent to 
 $f(x)\leq O(g(x))$.)
Here $a$ is any constant less than $a_0$, where the constants
$a_0=0.8168146...$ and $b=0.2587966...$ and $c=0.2064969...$ are defined as follows.
Let $B$ be the unique root of $1/B+\log B=1+\log 2$ in the interval $(0,1)$,
let $A$ be the unique root of $\sum_{n=1}^\infty A^n((n+1)\log(n+1)-n\log n-1)=1$ 
in the interval $(0,1)$, and put $a_0=-1/(2\log A)$ and $b=B\log 2$ and 
$c=(B\log 2)/(2-2B)$.
%
% A note on the existence of A and B: the function  theta(z) = 1/z + log(z)
% is decreasing on (0,1), since its derivative is  -1/z^2 + 1/z  which is negative,
% and  theta(1) = 1 < 1 + (log 2) < infinity = lim_{z->0+} theta(z).
% So there is indeed a unique z=B in (0,1) with  theta(B) = 1 + log 2.
%    We find  B = 0.3733646177016740842484484366...
% Next put  beta(z) = SUM_{n=1 to infinity} a_n z^n
% where  a_n = (n+1) log(n+1) - n log(n) - 1.
% Note that  a_n > 0  and  a_n is asymptotic to log(n);
% thus,  beta(z)  is defined for z in [0,1), and is strictly increasing,
% and  lim_{z->1-} beta(z) = infinity.  Since beta(0)=0, there is a unique z=A in (0,1) with
% beta(z)=1.
%    We find  A = 0.542598586098471021959...
%
The lower bound is proved by considering numbers $N$ which are products of
a fixed number $k$ of distinct primes, the least of which is $11$; for such $N$
we have $g_0(N)=1-2^{k-1}+(1/12)\prod_{p|N}(p+1)$, and the lower bound follows from
\cite{MP}.
The upper bound is proved by optimizing the choice
of $d$ in our density-zero proof, after one has modified that proof by
replacing the $o(x)$ bound from \cite[p.~356]{HW} by the more precise bound
from~\cite{HRa}.
Similar bounds have been proved for the number $\#V(x)$ of distinct values of Euler's
$\varphi$-function not exceeding $x$.  In this setting, our upper bound is essentially
an optimization of the paper \cite{Pi}.  In \cite{MP}, upper and lower bounds are proved
for $\#V(x)$, both of which have the same shape as the above lower bound (but the upper bound is
for any $b>b_0$).  The precise order of magnitude of $\#V(x)$ is determined in
the excellent paper \cite{Fo}, which (together with \cite{Fo2}) also contains several other
interesting results whose analogues would be interesting to study in our situation.  However, the 
multiplicativity of $\varphi(n)$ plays a crucial role in all proofs giving better upper 
bounds for $\#V(x)$ than our upper bound for $\#S(x)$ above.  Since $g_0(N)$ 
(and $g_0(N)-1$) is not multiplicative, it seems that these methods do not apply to 
$\#S(x)$.  We do not know which of our upper and lower bounds is closer to the truth.
\end{remark}

%####################################################################
%####################################################################

\section{Distribution of $g_0(N)$ modulo primes}
\label{sec-distribution}

We now study the distribution of $g_0(N)$ modulo a fixed
prime $\ell$ as $N$ varies.  
We will see that some residue classes occur more often
than others.  Let $\ell$ be a fixed prime.  We may restrict to levels 
$N$ having a prime factor $p$ congruent to $-1$ modulo $12\ell$, since 
the set of such $N$'s has density 1 in the set of positive integers 
(as shown by Dirichlet, cf.~\cite[p.~73]{Se}).
%
% Actually: in order to prove there are infinitely many primes congruent to -1 mod 12\ell,
% Dirichlet showed that  \sum_{primes p\equiv -1 \pmod {12\ell}} (1/p)  =  infinity.
% The desired density is clearly bigger than  1 - \prod_p (1 - 1/p)
% where the product is extended over any finite set of primes p congruent to -1 mod 12\ell.
% If follows that the desired density is 1.
%
Our assumption on $N$ forces $\nu_2=\nu_3=0$ and $12\ell\mid\mu$, so
\begin{equation*}
g_0(N)\equiv 1 - \nu_\infty/2\pmod{\ell}.
\end{equation*}
First consider the case $\ell=2$: for our restricted class of $N$'s, we 
have $g_0(N)\equiv 1\pmod{2}$ unless $N$ is either $p^{2r}$ or $4p^{2r}$,
so certainly $g_0(N)$ is odd with probability 1 (as was observed in 
Section~\ref{sec-data}).

We have seen that the sequence of residues mod 2 of $g_0(N)$ is biased.  We now show
a similar result modulo other primes $\ell$.  We compute the density of positive
integers $N$ for which $g_0(N)\equiv 1\pmod{\ell}$.  As above, we may assume
$N$ has a prime factor $p$ with $12\ell\mid (p+1)$; then $\ell\mid (g_0(N)-1)$
is equivalent to $\ell\mid\nu_\infty$, and this occurs precisely when either
$\ell^3\mid N$ or some prime congruent to $-1$~mod~$\ell$ divides $N$ with (positive)
even multiplicity.  By a standard argument (similar to \cite[Thm.~2.18]{Na}),
the probability that $N$ does not satisfy either of these conditions is
\begin{equation*}
 (1 - 1/\ell^3) \prod_{\substack{s\equiv-1\text{ (mod $\ell$)}\\ s \text{ prime} }} 
   (1-1/s^2+1/s^3-1/s^4+1/s^5-\dots),
\end{equation*}
which equals $(1 - 1/\ell^3) \prod (1 - 1/(s^2+s))$.  
The following table gives upper bounds for the probability $P(\ell)$ that 
$g_0(N)\equiv 1\pmod\ell$: 
\begin{center}
\begin{tabular}{c|c|c|c|c|c|c|c|c}
$\ell$&3&5&7&11&13&17&19&23 \\
\hline
$P(\ell)<$&1/4&1/78&1/105&1/653&1/1542&1/1793&1/978&1/5821
\end{tabular}
\end{center}
% Also P(29) < 1/10439.  And so on.
In these first few cases, we see that the probability is much less
than $1/\ell$.  More generally, for every $\ell$ we have 
$P(\ell)<3/\ell^2$ (where the constant `3' is not optimal):
\begin{align*}
P(\ell) &= 1 - (1 - 1/\ell^3)\prod_{\substack{s\equiv-1\text{ (mod $\ell$)}\\ 
s \text{ prime} }}  (1 - 1/(s^2+s))  \\
&< 1 - (1-1/\ell^3)\prod_{n=1}^\infty (1 - 1/((n\ell)^2-n\ell)) \\
&< 1 - (1-1/\ell^3)\prod_{n=1}^\infty (1 - 1.5/(n\ell)^2)  \\
&< 1 - (1-1/\ell^3)(1-\sum_{n=1}^\infty (1.5/\ell^2)/n^2) \\
% Here we use the fact that, if the a_i are positive and \sum a_i < 1,
% then \prod (1 - a_i) > 1 - \sum a_i.  This fact is easy to prove by
% induction.
&= 1 - (1-1/\ell^3)(1-\pi^2/(4\ell^2)) \\
%&< (\pi^2/4)\ell^{-2} + \ell^{-3} \\
&< 3\ell^{-2}.
\end{align*}

Also note that other cosets mod $\ell$ will have special behavior as well.
If $N$ is squarefree, $\nu_\infty$ is a power of 2, so that $g_0(N)$ is congruent 
modulo~$\ell$ to a number of the form $1-2^k$.  Recall that $N$ is squarefree 
with probability $6/\pi^2=0.6079...$~\cite[p.~269]{HW}.  Thus, if 2 is not a 
primitive root mod~$\ell$, then the residue classes (mod $\ell$) of the integers
$1-2^k$ will occur more frequently than the other residue classes.  For example,
for $\ell=7$, the classes $0,4,6$ occur much more frequently than the classes 
$2,3,5$.

%#######################################################################
%#######################################################################
%#######################################################################

% The following remark should only appear in electronic versions of this note:
% 
\vskip .3cm
\noindent
\scriptsize{Closing remark: readers wishing to see further details of various 
arguments in this note are invited to look at the comments in the \TeX\ source code.}

%#######################################################################
%#######################################################################
%#######################################################################

%\bibliographystyle{plain}

\end{document}